\documentclass[10pt]{article}
\usepackage{amsmath,amsthm,amsfonts,amssymb,amstext}
\theoremstyle{plain}

\newtheorem{theorem}{Theorem}[section]
\newtheorem{corollary}[theorem]{Corollary}
\newtheorem{lemma}[theorem]{Lemma}
\newtheorem{proposition}[theorem]{Proposition}
\newtheorem{definition}[theorem]{Definition}
\newtheorem{remark}[theorem]{Remark}

\newtheorem*{mainresult}{Main Result}
\newtheorem*{observation}{Observation}

\newtheorem*{question}{Question}

\begin{document}

\newcommand{\quasiOne}{\ensuremath{\pi^{G}_{H_{1}}}}
\newcommand{\quasiTwo}{\ensuremath{\pi^{G}_{H_{2}}}}

\newcommand{\EigenRep}{\ensuremath{\pi^{G}_{\lambda}}}
\newcommand{\EigenRepOne}{\ensuremath{\pi^{H_{1}}_{\lambda}}}
\newcommand{\EigenRepTwo}{\ensuremath{\pi^{H_{2}}_{\lambda}}}

\newcommand{\weak}{\ensuremath{\stackrel{w}{\sim}}}
\newcommand{\auto}{\ensuremath{\stackrel{a}{\sim}}}

\newcommand{\Z}{\mathbb Z}
\newcommand{\Q}{\mathbb Q}
\newcommand{\R}{\mathbb R}
\newcommand{\N}{\mathbb N}
\newcommand{\C}{\mathbb C}
\newcommand{\F}{\mathbb F}
\newcommand{\f}{\mathbb f}
\newcommand{\K}{\mathbb K}
\newcommand{\HYP}{\mathbb H}
\newcommand{\germ}{\mathfrak}
\newcommand{\g}{\mathfrak{g}}
\newcommand{\h}{\mathfrak{h}}
\newcommand{\W}{\mathcal{W}}
\newcommand{\T}{\mathcal{T}}
\newcommand{\D}{\mathcal{D}}
\newcommand{\A}{\mathcal{A}}

\newcommand{\Char}{\ensuremath{\operatorname{Char}}}
\newcommand{\Diff}{\ensuremath{\operatorname{Diff}}}
\newcommand{\Ind}{\ensuremath{\operatorname{Ind}}}
\newcommand{\Res}{\ensuremath{\operatorname{Res}}}
\newcommand{\DIM}{\ensuremath{\operatorname{dim}}}
\newcommand{\diag}{\ensuremath{\operatorname{diag}}}
\newcommand{\End}{\ensuremath{\operatorname{End}}}
\newcommand{\Isom}{\operatorname{Isom}}
\newcommand{\SU}{\operatorname{SU}}
\newcommand{\GL}{\operatorname{GL}}
\newcommand{\BC}{\operatorname{BC}}
\newcommand{\SO}{\operatorname{SO}}
\newcommand{\SL}{\operatorname{SL}}
\newcommand{\U}{\operatorname{U}}
\newcommand{\Tr}{\operatorname{Tr}}
\newcommand{\su}{\mathfrak{su}}

\newcommand{\Ad}{\operatorname{Ad}}
\newcommand{\ad}{\operatorname{ad}}
\newcommand{\Spec}{\operatorname{Spec}}
\newcommand{\spec}{\operatorname{spec}}
\newcommand{\Normal}{\operatorname{N}}
\newcommand{\Spin}{\operatorname{Spin}}
\newcommand{\Hol}{\operatorname{Hol}}
\newcommand{\Aut}{\operatorname{Aut}}
\newcommand{\Sim}{\operatorname{Sim}}


\title{Isospectral simply-connected homogeneous spaces \linebreak and the spectral rigidity of group
actions}
\author{Craig J. Sutton}
\date{}
\maketitle

\begin{abstract} 
\noindent
We generalize Sunada's method to produce new examples of closed, locally non-isometric manifolds 
which are isospectral. In particular, we produce pairs of isospectral, simply-connected, 
locally non-isometric normal homogeneous spaces. These pairs also allow us to see that in general group
actions with discrete spectra are not determined up to measurable conjugacy by their spectra. 
In particular, we show this for lattice actions. 
\end{abstract}



\tableofcontents


\section{Introduction}

Spectral geometry is the study of the relationship between the geometry of a Riemannian manifold $(X,m)$ 
and the spectrum of the associated Laplace-Beltrami operator $\Delta$ acting on $C^{\infty}(X)$.
Specifically, one is concerned with the extent to which the spectrum encodes geometric information.
While the spectrum does determine some geometric properties such as total scalar curvature, volume, 
and dimension, in general it does not determine a Riemannian manifold up to isometry. This was demonstrated
for the first time by Milnor in 1964 when he produced examples of 16-dimensional tori which are isospectral
yet non-isometric \cite{Milnor}.  Hence, in order to better understand the interplay between the geometry of a
Riemannian manifold and its spectrum, other such examples must be studied.

During the past two decades many new non-isometric isospectral spaces have been found 
(e.g., \cite{GW}, \cite{BT}, \cite{BG}, \cite{Gornet1}, \cite{Gornet2}, \cite{Gordon}, \cite{Szabo} and
\cite{GGSWW}).\footnote{For a more comprehensive discussion of the spectral geometry landscape the
reader is encouraged to see \cite{Gordon3}.}   
The first examples of topological significance were produced by Vign\'{e}ras and Ikeda. 
In \cite{Vigneras} examples of 3-dimensional hyperbolic spaces with non-isomorphic fundamental groups were
constructed and in \cite{Ikeda} isospectral lens spaces were produced. These examples demonstrated for the
first time that the topology of the manifold is not a spectral invariant. However, it is worth noting
that these isospectral spaces (along with all other understood examples) have a common universal cover.

Inspired by a result from number theory, Sunada produced the first general method for constructing 
pairs of isospectral manifolds.

\begin{theorem}[Sunada's Method, \cite{Sunada}]
Let $(X,m)$ and $(X_{0}, m_{0})$ be Riemannian manifolds and $\pi : X \to X_{0}$ be a finite Riemannian covering 
with covering transformation group $G$.
Now suppose $\pi_{1} : X_{1} \to X_{0}$ and $\pi_{2} : X_{2} \to X_{0}$ are the Riemannian coverings corresponding 
to subgroups $\Gamma_{1}, \Gamma_{2} \leq G$ respectively. If for every $g\in G$ we have 
$\#(g^{G} \cap \Gamma_{1}) =\#(g^{G} \cap \Gamma_{2})$ (where $g^{G}$ denotes the conjugacy class of
$g$),  then $X_{1}$ and $X_{2}$ are isospectral.
\end{theorem}

Many of the examples of isospectral pairs that arise in the literature can be explained by Sunada's method or 
one of its generalizations. One generalization that will be of interest to us is the following. 

\begin{theorem}[Sunada-Pesce Method, \cite{Pesce}] 
Let $(X,m)$ be a Riemannian manifold, $G \leq \Isom (X,m)$ 
closed, $K$ the generic stabilizer of the action of $G$ on $X$ (see Section~\ref{Sec:TheMethod}), 
and $\Gamma_{1}, \Gamma_{2} \leq G$ be discrete such 
that the manifolds $\Gamma_{1} \backslash X$ and $\Gamma_{2} \backslash X$ are compact. If the quasi-regular representations 
$\pi_{\Gamma_{1}}^{G}$ and $\pi_{\Gamma_{2}}^{G}$ of $G$ are $K$-equivalent (see Section~\ref{Sec:RelEquivalence}), then the 
Riemannian quotients $(\Gamma_{1} \backslash X , m_{1})$ and $(\Gamma_{2} \backslash X , m_{2})$ are isospectral.
\end{theorem}
\noindent
As with all previous generalizations of Sunada's method pairs arising in this manner are not simply-connected.
Also, the resulting pairs have a common Riemannian covering, namely $X$, and consequently are locally isometric.
This causes us to wonder whether one can generalize Sunada's method so that it produces locally
non-isometric, simply-connected isospectral pairs.

A natural approach to this would be to take quotients of simply-connected Lie groups by non-trivial 
connected subgroups, which leads us to the following long standing problem in the spectral geometry
community.

\begin{question}\label{Ques:LongStanding}
Are there examples of Riemannian manifolds $(X,m)$ such that one can find $H_{1}, \; H_{2} \leq \Isom(X,m)$ non-trivial 
and connected such that the quotient manifolds $(X/H_{1}, m_{1})$ and $(X/H_{2}, m_{2})$ are isospectral yet non-isometric?
\end{question}

In this paper we are able to answer this question positively. Using a generalization of the 
Sunada-Pesce method and a result of Larsen and Pink \cite{LP} we establish the following.

\begin{mainresult}[Theorem \ref{Thm:MyExamples}]
There exists a connected, simply-connected semi-simple real Lie group $H$ which for infinitely 
many $n \in \N$ admits reducible faithful representations 
$\rho_{1}, \; \rho_{2} : H \to \SU(n)$, where $\rho_{1} \not\auto \rho_{2}$
(see Definition~\ref{autoequivalence}) and 
$H_{1} = \rho_{1}(H)$ and $H_{2} = \rho_{2}(H)$ are not conjugate by $\Aut(\SU(n))$.
If we equip $\SU(n)$ with a bi-invariant metric $m$, then the simply-connected, normal homogeneous spaces 
$(\SU(n)/H_{1}, m_{1})$ and $(\SU(n)/H_{2}, m_{2})$ are isospectral yet locally non-isometric.
\end{mainresult} 

Recently, we have learned that Schueth has also obtained examples of isospectral homogeneous 
spaces \cite{Schueth2}. In fact, she produces a continuous family of pairwise isospectral left-invariant
metrics on a simply-connected Lie group. Our examples can be distinguished from Schueth's in that they are
normal homogeneous spaces; that is, they have the metric  induced by the bi-invariant metric on $G$, 
and our spaces are quotients of
$G$ by \emph{non-trivial} connected subgroups $H_{1}, H_{2} \leq G$, which are representation  equivalent
(see Definition~\ref{Def:Equiv}).
They can also be distinguished from Schueth's examples in that the method
of construction necessitates an enormous dimension for the resulting homogeneous spaces. An
estimate shows the simplest example to have a dimension on the order of $10^{10}$ .

The spaces constructed in this paper along with those of Schueth \cite{Schueth, Schueth2}, Gordon
\cite{Gordon2} and Szab\'{o} \cite{Szabo2} are the only known examples of closed, \emph{simply-connected}, 
locally non-isometric isospectral spaces. Schueth, Gordon and Szab\'{o} construct their isospectral spaces by
fixing a particular simply-connected manifold and then creating isospectral metrics on this space through
various interesting techniques. Consequently, the resulting isospectral spaces are always homeomorphic.
At the present time it is unclear to the author whether the isospectral pairs presented in this paper are
homeomorphic. A negative answer would demonstrate for the first time that the universal cover is not a spectral 
invariant.\footnote{In looking for simpler candidates for non-homeomorphic isospectral simply-connected spaces one
might consider the Aloff-Wallach spaces \cite{AW} as \emph{normal} homogeneous spaces. These spaces are
simply-connected, however, it can be shown that isospectral normal homogeneous Aloff-Wallach spaces are
necessarily isometric and hence homeomorphic (see \cite{Blanchard}, \cite{JLPR} or \cite{Sutton}).} 

The structure of this paper is as follows. In Section~\ref{Sec:SunadaPesceMethod} we will discuss 
the proof of the generalized Sunada-Pesce method. We will use this method in Section~\ref{Sec:NewExamples}
to construct new examples of isospectral yet locally non-isometric pairs of Riemannian manifolds. In
Section~\ref{Sec:FiberBundles} we will establish a method for constructing isospectral fiber bundles 
with isospectral fibers. Finally, we recall that a well-known theorem of von Neumann states that
two actions of an abelian locally compact group with \emph{discrete} spectra are measurably conjugate
if the actions are isospectral \cite{vonNeumann}.  In Section~\ref{Sec:GroupActions} we will use the examples
constructed in  Section~\ref{Sec:NewExamples} to demonstrate that group actions with discrete spectra
are not classified up to measurable conjugacy by their spectra and hence von Neumann's result is not true in
general.

\medskip
\noindent
{\bf Notation.} We will use the following notation.
\begin{enumerate}
\item We will let $\Spec(\Delta)$ denote the spectrum of the Laplacian taking multiplicities into account.
\item Given a representation $\rho : G \to \GL(V)$ of a Lie group $G$ we will let 
$\Res^{G}_{H}(\rho)$ denote the restriction of $\rho$ to $H$ for any subgroup $H \leq G$.
\item Given a Lie group $G$, a closed subgroup $H \leq G$ and a representation $\tau : H \to \GL(V)$ we will let 
$\Ind^{G}_{H}(\tau)$ denote the induced representation (see p. \pageref{induced} ).
\item Given two representations $\rho_{1}: G_{1} \to \GL(V_{1})$ and 
$\rho_{2}: G_{2} \to \GL(V_{2})$ we let 
$\rho_{1} \otimes \rho_{2} : G_{1} \times G_{2} \to \GL(V_{1} \otimes V_{2})$ denote the {\bf outer tensor
product} given by $\rho_{1} \otimes \rho_{2}(g_{1},g_{2}) = \rho_{1}(g_{1}) \otimes \rho_{2}(g_{2})$.
\item Given two representations $(\rho, V)$ and $(\tau, W)$ of a Lie group $G$ we will let 
$[\rho: \tau]$ denote the multiplicity of $\tau$ in $\rho$. In the case where $E$ is a field extension of $F$
we will let $[E:F]$ denote the degree of the extension.
\item We will use the symbol ``$\leq$'' to denote both vector subspaces and subgroups.
\end{enumerate}

\noindent
{\bf Acknowledgements.} The work presented in this paper is part of my thesis \cite{Sutton} carried out at the 
University of Michigan. It is my great pleasure to thank my advisor, Ralf Spatzier, for introducing me to the area
of spectral geometry and, more importantly, for being a generous and supportive mentor. I am also indebted to
Gopal Prasad for making me aware of the work of Larsen and Pink concerning dimension data and to Krishnan Shankar
for discussing homogeneous spaces with me. Thanks also go to the referee for providing helpful comments
concerning the exposition of this article.
\smallskip


\section{Generalized Sunada-Pesce Method}\label{Sec:SunadaPesceMethod}

Developing techniques for constructing isospectral manifolds is one of the central concerns of inverse spectral geometry. 
The examples these techniques yield allow us to discover the geometric data that cannot be recovered from
the spectrum of the  Laplacian. In this section we will generalize Sunada's well-known method for
constructing pairs of isospectral manifolds. More specifically, we will generalize the Sunada-Pesce method
to allow one to obtain isospectral pairs by taking  quotients by non-trivial connected groups. 
By considering such quotients we open up the possibility that the resulting isospectral pairs need not have a
common Riemannian covering or common universal cover, which is impossible under other 
versions of Sunada's method.
In fact, in Section~\ref{Sec:NewExamples} we will show that through
this method we can construct many pairs of isospectral, simply-connected, locally non-isometric spaces. 
We begin by reviewing the concept of relatively equivalent representations.

 
\subsection{Relative Equivalence}\label{Sec:RelEquivalence}

Two representations $\rho_{1}: G \to \GL(V_{1})$ and $\rho_{2}: G \to \GL(V_{2})$ of a Lie
group $G$ are said to be {\bf equivalent}, denoted by $\rho_{1} \sim \rho_{2}$ or $(\rho_{1}, V_{1}) \sim
(\rho_{2}, V_{2})$, if there exists a vector space isomorphism
$T: V_{1} \to V_{2}$ such that $\rho_{2}(g) \circ T = T \circ \rho_{1}(g)$ for any $g \in G$.
Now, consider a representation $\rho : G \to \GL(V)$ of a unimodular Lie group $G$. 
For any compact subgroup $K$ of $G$ and any representation $(\tau, W)$ of $K$ we may consider the
vector subspace $V^{\tau} \equiv \oplus V_{\alpha} \leq V$, where the direct sum is taken over $K$-invariant
subspaces $V_{\alpha} \leq V$  such that $(\Res^{G}_{K}(\rho),V_{\alpha})\sim (\tau, W)$. The representation 
$(\Res^{G}_{K}(\rho), V^{\tau})$ can then be extended to a representation  of $G$ by considering the vector space
$$V_{\tau} \equiv \cap \{ L \leq V : V^{\tau} \leq L \mbox{ and } L \mbox{ is $G$-invariant} \}.$$
This subrepresentation will be denoted by $(\rho_{\tau}, V_{\tau})$. For any two representations $(\rho, V)$ 
and $(\tilde{\rho}, \widetilde{V})$ 
of $G$ we will agree to say they are \textbf{$\tau$-equivalent}, denoted by $\rho \sim_{\tau} \tilde{\rho}$ or 
$(\rho, V) \sim_{\tau} (\tilde{\rho},\widetilde{V})$,  
if the subrepresentations $(\rho_{\tau}, V_{\tau})$ and $(\tilde{\rho}_{\tau},\widetilde{V}_{\tau})$ are
equivalent.

In this paper we will be concerned with representations of $G$ which are $1_{K}$-equivalent, 
where $1_{K}$ denotes the trivial representation of $K$ on $\C$. To be consistent with \cite{Pesce} we will refer to
this as \textbf{$K$-equivalence}, and we will  denote $V^{1_{K}}$ and $(\rho_{1_{K}}, V_{1_{K}})$ by $V^{K}$ and
$(\rho_{K}, V_{K})$ respectively. As usual we will let 
$\widehat{G}$ denote the set of equivalence classes of irreducible representations of $G$ and we will agree to let 
$\widehat{G}_{K}$ denote the set of equivalence classes of representations of $G$ which admit 
non-trivial $K$-fixed vectors; that is, 
$\widehat{G}_{K} \equiv \{ [(\rho, V)] \in \widehat{G} : V_{K} \neq 0 \}$.   
We will illustrate the concept of $K$-equivalence by constructing two representations of a group $G$,  which are
$K$-equivalent for some subgroup $K \leq G$, but not equivalent. In doing this it will be useful to review the
method of induction.

Let $G$ be a locally compact group, $H$ a closed subgroup and $\rho: H \to \GL(V)$ a unitary representation of $H$.
The representation $\rho$ gives a representation of $G$ known as the {\bf induced representation}\label{induced},
denoted by
$\Ind^{G}_{H}(\rho)$, which acts on the vector space
$$ \widetilde{V} \equiv \{f: G \rightarrow V \mbox{ $L^{2}$-function} : f(gh)= \rho(h)^{-1}f(g)
\mbox{ for all } h \in H , \mbox{ almost every } g \in G \}$$
by 
$$(\Ind^{G}_{H}(\rho)(g)f)(x) \equiv f(g^{-1}x)$$
for any $g\in G$ and $f \in \widetilde{V}$.
In the case where $\rho = 1_{H}$ is the trivial representation of $H$, $\Ind^{G}_{H}(\rho)$ is 
a representation of $G$ on $L^{2}(G/H)$ known as the \textbf{quasi-regular representation} of $G$ 
with respect to $H$, which we will denote by $\pi^{G}_{H}$. We now make the following observation.

\begin{observation}\label{Obs:Kequivalence}
Let $G$ be a Lie group. Now consider subgroups $K, H_{1}, H_{2} \leq G$ (closed) such that 
$K\backslash G / H_{i} = \{\mbox{point}\}$ (equivalently $H_{i}$ acts transitively on $K\backslash G$) 
for $i = 1,2$, then $\pi_{H_{1}}^{G} \sim_{K} \pi_{H_{2}}^{G}$. If $G$ is compact and 
$\mbox{dim}(G/H_{1}) \neq \mbox{dim}(G/H_{2})$, then we may conclude $\pi_{H_{1}}^{G} \not\sim \pi_{H_{2}}^{G}$.
\end{observation}

For $n \geq 2$, we can see that if $G = \SO(4n), \; H_{1} = U(2n), \; H_{2} = Sp(n)$ and $K = \SO (4n-1)$, then the above 
implies that the representations $\pi_{H_{1}}^{G}$ and $\pi_{H_{2}}^{G}$ are $K$-equivalent, yet inequivalent.
Indeed, we view $K$ as a subgroup of $G$ under the imbedding $A \mapsto [1] \oplus A$ and we consider $H_{1}$ and $H_{2}$ as subgroups 
of $G$ by using the standard imbedding of complex and quaternionic matrices into the real matrices (see \cite[p. 34-36]{Knapp1}). 
Then since $H_{1}$ and $H_{2}$ act transitively on $S^{4n-1} = K \backslash G$ we see that 
$K\backslash G / H_{i} = \{\mbox{point} \}$ for $i=1,2$. Hence, $\pi_{H_{1}}^{G} \sim_{K} \pi_{H_{2}}^{G}$, however, 
since $\mbox{dim}(G/H_{1}) \neq \mbox{dim}(G/H_{2})$ (for $n\geq 2$) we find that 
$\pi_{H_{1}}^{G} \not\sim \pi_{H_{2}}^{G}$.

We conclude this section with a little jargon.

\begin{definition}\label{Def:Equiv}
Let $G$ be a compact Lie group and $K\leq G$ compact. We will say that two closed subgroups $H_{1}, H_{2} \leq G$ are
\begin{enumerate}
\item {\bf Representation equivalent} if $\pi_{H_{1}}^{G} \sim \pi_{H_{2}}^{G}$.
\item {\bf K-equivalent} if $\pi_{H_{1}}^{G} \sim_{K} \pi_{H_{2}}^{G}$.
\end{enumerate}
\end{definition}

 
\subsection{The Method}\label{Sec:TheMethod}

Before stating our method for constructing isospectral Riemannian manifolds, we recall 
the notion of the generic stabilizer. 

\begin{definition}\label{Stabilizer}
Suppose $G$ is a Lie group which has a proper $C^{\infty}$-action on a manifold $X$. For each 
$x\in X$ let $G_{x}$ denote the stabilizer of $x$. There exists a subgroup $K$ of $G$ called the 
\textbf{generic stabilizer} with the following properties:
\begin{enumerate}
\item For all $x \in X$, $K$ is conjugate to a subgroup of $G_{x}$.
\item There exists an open and dense subset $U$ in $X$  such that for all $x \in U$ $K$ and $G_{x}$ are conjugate.
\end{enumerate}
Orbit spaces of the type $G/K$ are known as \textbf{principal orbits}. 
\end{definition}

With this terminology we may now state the following proposition.

\begin{theorem}[Generalized Sunada-Pesce Technique] \label{Thm:GenSuPe}
Let $(X, m)$ be a compact Riemannian manifold and $G \leq \Isom(X,m)$ a compact Lie 
group. We will let $K$ denote the generic stabilizer of the action of $G$ on $X$. Now suppose 
$H_{1}, H_{2} \leq G$ are closed, $K$-equivalent subgroups which act freely on $X$ and are such that 
the Riemannian submersions $\pi_{i}:X \rightarrow X/H_{i}$, $i=1,2$, have totally geodesic fibers. 
It then follows that the Riemannian  manifolds $X/H_{1}$ and $X/H_{2}$ are isospectral on functions.
\end{theorem}

\noindent
As in Pesce's original paper \cite{Pesce} the proof of this theorem is an application of Frobenius' 
reciprocity theorem and a result of Donnelly. But first we recall the following result.

\begin{lemma}\label{Lemma:TotallyGeodesicFibers}
Let $(E^{k},m)$ and $(B^{n}, m_{B})$ be Riemannian manifolds.
Let $\pi : (E, m) \to (B, m_{B})$ be a Riemannian submersion with totally geodesic fibers, then the 
eigenfunctions of $B$ are functions whose pullbacks are eigenfunctions on $E$. In fact, if $f$ is an
eigenfunction  of $\Delta_{B}$ with eigenvalue $\lambda$, then its pullback $f\circ \pi$ is an eigenvalue
of $\Delta_{E}$  with eigenvalue $\lambda$. Hence, we see that $\Spec(\Delta_{B}) \subset
\Spec(\Delta_{E})$. 
\end{lemma}

\begin{proof}
Let $f \in L^{2}(B)$ be such that $\Delta_{B} f = \lambda f$ and let $\tilde{f} = f \circ \pi \in L^{2}(E)$ be its
pullback to $E$. 
Now fix $x \in E$ and let $\{ \gamma_{1}, ..., \gamma_{k}\}$  be a collection of geodesics such that $\gamma_{i} (0) =
x$ for all
$i$ and 
$\{\dot{\gamma}_{1}(0),...,\dot{\gamma}_{k}(0) \}$ is an orthonormal basis for $T_{x}E$ with 
$\{\dot{\gamma}_{1}(0),...,\dot{\gamma}_{n}(0) \}$ vertical (that is, tangent to the fiber through $x$).
Then 
\begin{equation}
\begin{split}
\Delta_{E} \tilde{f} (x) &= - \sum ^{k}_{1} \frac{d^2}{dt^{2}} (\tilde{f}\circ \gamma_{i})(0) \\
			 &= - \sum ^{k}_{n+1} \frac{d^2}{dt^{2}} (\tilde{f}\circ \gamma_{i})(0) \\
			&= - \sum ^{k}_{n+1} \frac{d^2}{dt^{2}} (f\circ ( \pi \circ \gamma_{i})) (0) \\
			&= \Delta_{B} f (\pi (x)) \\
			&= \lambda \tilde{f} (x). \\
\end{split}
\end{equation}
This shows us that pullbacks of eigenfunctions on $B$ are eigenfunctions on $E$ with the same eigenvalue.
So, we obtain $\Spec (\Delta_{B}) \subset \Spec (\Delta_{E})$.
\end{proof}

\begin{proof}[Proof of Theorem~\ref{Thm:GenSuPe}]
Let $\Delta , \Delta_{1},$ and $\Delta_{2}$ denote the Laplace-Beltrami operator on $X, X/H_{1}$, and $X/H_{2}$ 
respectively. Since $\pi_{i} : X \to X/H_{i}$ ($i = 1,2$) has totally geodesic fibers it follows 
from Lemma~\ref{Lemma:TotallyGeodesicFibers} that 
$\Spec (\Delta_{i}) \subset \Spec (\Delta)$ for $i=1,2$. We also recall that the action of $\Isom (X,m)$ on 
$L^{2}(X)$ commutes with the Laplacian. Hence, $\Isom (X,m)$ preserves the eigenspace decomposition of 
$L^{2}(X)$. So for any $\lambda \in \Spec (\Delta) \mbox{ and } H \leq \Isom ( X, m)$ (closed) we have a 
representation $\pi^{H}_{\lambda}$ of $H$ on $L^{2}(X,m)_{\lambda}$ given by 
 $\pi^{H}_{\lambda} (h) .f = f \circ h^{-1}$. In our situation we will be interested in 
$\EigenRep, \EigenRepOne, \mbox{ and } \EigenRepTwo$ for $\lambda \in \Spec (\Delta)$.

Now it is clear that for any  $\lambda \in \Spec (\Delta)$ we have 
$\dim L^{2}(X/H_{i})_{\lambda} = [\pi^{H_{i}}_{\lambda} : 1_{H_{i}}]$. 
Indeed, for $H \leq \Isom(X,m)$ (closed) we let $L^{2}(X)^{H} = \{ f \in L^{2}(X) : h.f = f \mbox{ for all } h \in H \}$.
One can see that $L^{2}(X/H) = L^{2}(X)^{H}$ and it follows that 
$L^{2}(X/H)_{\lambda} = L^{2}(X)^{H}_{\lambda}$. Hence, $\DIM L^{2}(X/H)_{\lambda} = [\pi^{H}_{\lambda} : 1_{H}]$.
We may now conclude that $(X/H_{1}, m_{1}) \mbox{ and } (X/H_{2}, m_{2})$ are isospectral if and only if 
$[\pi^{H_{1}}_{\lambda} : 1_{H_{1}}] = [\pi^{H_{2}}_{\lambda} : 1_{H_{2}}]$ for all $\lambda \in \Spec (\Delta ) $.
Since it is clear that for every $H \leq G$ $\pi^{H}_{\lambda} = \Res^{G}_{H}(\EigenRep)$, Frobenius 
reciprocity gives us the following:

\begin{equation}
\begin{split}
[\pi^{H_{i}}_{\lambda} : 1_{H_{i}}]	&= [\Res^{G}_{H_{i}}(\EigenRep) : 1_{H_{i}}] \\
					&= [\Res^{G}_{H_{i}}(\sum_{\rho \in \widehat{G}}[\EigenRep:\rho]\rho) : 1_{H_{i}}] \\
					&= \sum_{\rho \in \widehat{G}}[\EigenRep:\rho][\Res^{G}_{H_{i}}(\rho) : 1_{H_{i}}] \\
					&= \sum_{\rho \in \widehat{G}}[\EigenRep:\rho][\Ind^{G}_{H_{i}}(1_{H_{i}}) : \rho] \\
					&= \sum_{\rho \in \widehat{G}}[\EigenRep:\rho][\pi^{G}_{H_{i}} : \rho]. \\
\end{split}
\end{equation}
We now recall the following theorem of Donnelly.

\begin{theorem}[\cite{Donn}, p. 25]
Let $G$ be a compact Lie group and $X$ a compact, smooth $G$-space with principal orbit type $G/K$; that is, $K$ is the 
generic stabilizer of the $G$-action on $X$. Then the decomposition of $L^{2}(X)$ into 
$G$-irreducibles contains precisely those finite dimensional representations appearing in the decomposition of 
$\pi^{G}_{K} = \Ind^{G}_{K}(1_{K})$ the quasi-regular representation of $G$ with respect to $K$. Also, if the orbit 
space $X/G$ has dimension greater than 1, then each irreducible appears an infinite number of times.
\end{theorem}

Now by Frobenius reciprocity we have $[\pi^{G}_{K}: \rho] = [ \Res^{G}_{K}(\rho) : 1_{K}]$ for each $\rho \in \widehat{G}$.
So we conclude from the above theorem that $[\pi^{G}_{\lambda}: \rho] \neq 0$ for some 
$\lambda \in \Spec (\Delta )$ if and only if $\rho_{K}$ is non-trivial. Consequently, for $i=1,2$, we see that 
$$[\pi^{H_{i}}_{\lambda} : 1_{H_{i}}] = \sum_{\rho \in \widehat{G}_{K}}[\EigenRep:\rho][\pi^{G}_{H_{i}} : \rho].$$ 
Finally, as a result of the $K$-equivalence of $\quasiOne$ and $\quasiTwo$ we obtain isospectrality.
\end{proof}

\section{Building New Examples}\label{Sec:NewExamples}

In this section we will use Theorem~\ref{Thm:GenSuPe} along with a result of Larsen and 
Pink \cite{LP} to produce the first pairs of non-isometric isospectral manifolds which are of the form 
$(X/H_{1}, m_{1})$ and $(X/H_{2}, m_{2})$, where $H_{1}, H_{2} \leq \Isom(X,m)$ are nontrivial and connected.
In particular, we will obtain the first examples of isospectral simply-connected, locally non-isometric, normal 
homogeneous spaces. 

\medskip

We begin by introducing a slightly more general notion of equivalence of representations.

\begin{definition}\label{autoequivalence}
Two representations $\tau_{1} : G \to \GL(V_{1})$ and $\tau_{2} : G \to \GL(V_{2})$ are said to be {\bf
automorphically equivalent}, denoted 
$\tau_{1} \auto \tau_{2}$ or $(\tau_{1}, V_{1}) \auto (\tau_{2}, V_{2})$, if there
exists a Lie group automorphism $\alpha : G \to G$ and a vector space isomorphism $T : V_{1} \to
V_{2}$  such that $T \circ \tau_{1}(g) = \tau_{2}(\alpha(g))\circ T$ for all $g \in G$. 
\end{definition}
\noindent
Clearly, equivalence implies automorphic equivalence just by taking $\alpha$ to be the identity map. 
However, a dramatic difference between these two definitions can be obtained by considering the 
irreducible representations of the additive group $\R$. 
For each $\theta \in \C$ we obtain an irreducible representation of $\R$ on $\C$ given by
$\pi_{\theta}(x)v = e^{2\pi i \theta x}v$ for any $x \in \R$ and $v \in \C$.
These are all the inequivalent irreducible representations of $\R$, but we see that for $\theta, \vartheta \in \R
\backslash
\{0\}$ $\pi_{\theta}(x) = \pi_{\vartheta}(\frac{\theta}{\vartheta}x)$, hence all of the non-trivial
irreducible representations are automorphically equivalent.\footnote{This example was pointed out to the
author by A. Knapp. A less dramatic example is obtained by comparing the standard representation of 
$i: \SU(n) \hookrightarrow \GL_{n}(\C)$ and $\sigma : \SU(n) \to \GL_{n}(\C)$ given by 
$\sigma(g) = \overline{i(g)}$, where the bar denotes complex conjugation.}

Now consider $G$ a connected, complex reductive Lie group and 
let $\rho: G \to \GL(V)$ be a faithful representation of dimension $n$. The \textbf{dimension data} of 
$(\rho, V)$ is defined as 
$$\{ (\sigma : GL(V) \rightarrow GL(W), \dim W^{G}): \sigma \mbox{ is a homomorphism and} \DIM W < \infty \}.$$
The objective of \cite{LP} is to determine the extent to which the \textbf{dimension data} of $(\rho,
V)$, determines the group $G$ and/or the representation $\rho : G \to \GL(V)$. The main result of their paper is
the  following.

\begin{theorem}[ \cite{LP}, p. 377] \label{Thm:LarsenPink}
Let $G$ be a connected, complex Lie group with $(\rho, V)$ a finite dimensional faithful representation. 
Then
\begin{enumerate}
\item The dimension data determine $G$ up to isomorphism. That is, if $\tau: G' \to \GL(V)$
is another representation of a connected, complex Lie group $G'$ with the same dimension data, 
then $G$ and $G'$ are isomorphic as Lie groups. (Notice that $G$ and $G'$ act on the same vector space
$V$.)
\item If $\rho$ is irreducible, the dimension data determine $\rho$ up to automorphic equivalence. That is, if
there  exists another faithful irreducible representation $\tau: G \to \GL(V)$ of $G$ with the same dimension data, 
then $\rho$ and $\tau$ are automorphically equivalent.
\item There exists a $G$ as above which admits a countably infinite number of pairs of reducible
representations
$(\rho_{1}, V)$ and $(\rho_{2}, V)$ of $G$, where $V$ is of arbitrarily large dimension, such that $\rho_{1}$ and 
$\rho_{2}$ have the same dimension data and $\rho_{1} \not\auto \rho_{2}$.
\end{enumerate}
\end{theorem}

\begin{remark}
In \cite{LP} the term automorphically equivalent is not used. Instead they use isomorphic. 
We have introduced this term so as not to cause confusion with the usual notion of equivalence, which is
sometimes referred to as isomorphic. 
\end{remark}

Our interest lies in the third part of the above theorem. 
We note that the method Larsen and Pink employ to produce the automorphically inequivalent pairs of
representations  with the same dimension data actually yields self-dual representations. That is, in the above 
$\rho_{1} \sim \rho^{*}_{1}$ and $\rho_{2} \sim \rho^{*}_{2}$, 
where for any representation $\tau : G \to \GL(W)$ $\tau^{*}(g) \equiv \tau(g^{-1})^{t}$ is the {\bf contragredient}
representation. Indeed, we recall that on p. 392 of \cite{LP} the group $G$ is constructed as a product of
non-isomorphic semisimple Lie groups 
$G_{1}, \ldots , G_{r}$ whose root systems $\Phi_{1}, \ldots , \Phi_{r} \subset \BC_{n}$ are subsystems of maximal 
rank $n$. They then choose formal characters 
$v_{1}, \ldots , v_{r} \in \Z[\Lambda_{\BC_{n}}]^{\W_{n}} \equiv \Z[\Z^{n}]^{\W_{n}}$, 
where $\W_{n} = \{\pm 1\}^{n} \rtimes S_{n}$ is the Weyl group of $\BC_{n}$ and $S_{n}$ is the permutation group on $n$
elements, such that for all $i, j = 1, \ldots , r$  there exists a faithful representation 
$\rho_{ij} : G_{i} \to \GL(V_{ij})$ with formal character $v_{j}$. Since the formal character $v_{j}$ is 
invariant under $\{\pm 1\}^{n} \rtimes S_{n}$ it follows that if $\lambda$ is a weight of 
$v_{j}$, then so is $-\lambda$. Hence, any representation with formal character $v_{j}$ is self-dual.
Larsen and Pink then consider the faithful representations 
$\rho_{1} = \oplus_{\sigma \in A_{r}} \rho_{1\sigma(1)} \otimes \cdots \otimes \rho_{r\sigma(r)} :
G \to \GL(V_{1})$ and 
$\rho_{2} = \oplus_{\sigma \in S_{r} - A_{r}} \rho_{1\sigma(1)} \otimes \cdots \otimes
\rho_{r\sigma(r)} : G \to \GL(V_{2})$,  where 
$$V_{1}  = \bigoplus_{\sigma \in A_{r}} V_{1\sigma(1)}\otimes \dots \otimes V_{r\sigma(r)}$$
and 
$$V_{2}  = \bigoplus_{\sigma \in S_{r}-A_{r}} V_{1\sigma(1)}\otimes \dots \otimes V_{r\sigma(r)}.$$
It is then clear that $\rho_{1} \sim \rho^{*}_{1}$ and $\rho_{2} \sim \rho^{*}_{2}$ and that 
$V_{1} \approx V_{2} \equiv V$. The representations $(\rho_{1}, V)$ and $(\rho_{2}, V)$ are the 
representations alluded to in Theorem~\ref{Thm:LarsenPink}(3).

If one now considers compact real forms we see that part three of Theorem~\ref{Thm:LarsenPink} can be 
recast as follows.

\begin{corollary}\label{Cor:Reduction}
There exists a compact, connected, semisimple real Lie group $H$ such that for infinitely many $n \in \N$ there 
exist faithful representations $\rho_{1},\; \rho_{2}: H \to \SU (n)$ with the same dimension data and such that 
$\rho_{1} \not\auto \rho_{2}$ and $\rho_{1} \not\auto \rho^{*}_{2}$. In fact, $H_{1} = \rho_{1}(H)$ and $H_{2} =
\rho_{2}(H)$  are not conjugate by $\Aut (\SU(n))$; that is, there are no automorphisms $\alpha$ of $\SU(n)$ such
that 
$\alpha(H_{1}) = H_{2}$.
\end{corollary}

\begin{proof}
The first part of this theorem is standard representation theory and follows for example from \cite[Theorem 4.11.14]{Var}.
As for the statement concerning the non-conjugacy of $H_{1}$ and $H_{2}$ we recall the following.

\begin{proposition}[see p. 56 of \cite{On2}]\label{Prop:Conjugacy}
Let $G$ be a connected, simple, non-abelian compact Lie group and $H$ a connected and simply-connected Lie group.
\begin{enumerate}
\item Let $\sigma, \; \tau : H \to G$ be two homomorphisms with discrete kernels. Then there exists 
$\alpha \in \Aut(G)$ such that 
$\alpha( \sigma(H)) = \tau(H)$ if and only if $\tau = \alpha \circ \sigma \circ \beta$ for a  certain 
$\beta \in \Aut(H)$.
\item Two homomorphisms $\tau, \sigma : H \to \SU(n)$ are {\bf conjugate by $\Aut(\SU(n))$} if and only if 
$\tau \sim \sigma$ or $\tau \sim \sigma^{*}$. Here conjugate by $\Aut(\SU(n))$ means 
there exists $\alpha \in \Aut(\SU(n))$ such that $\tau = \alpha \circ \sigma$.
\end{enumerate}
\end{proposition}

Now let's suppose $\alpha(H_{1}) = H_{2}$ for some $\alpha \in \Aut (\SU(n))$.
Then the first part of the above shows us that $\rho_{2} = \alpha \circ \rho_{1} \circ \beta$ for some 
$\beta \in \Aut(H)$. The second part of Proposition~\ref{Prop:Conjugacy} then implies that 
$\rho_{2} \sim \rho_{1} \circ \beta$ or $\rho_{2} \sim \rho^{*}_{1} \circ \beta$. 
That is, $\rho_{2} \auto \rho_{1}$ or $\rho_{2} \auto \rho^{*}_{1}$,  which is a contradiction. 
Hence, $H_{1}$ and $H_{2}$ are not conjugate by $\Aut(\SU(n))$.
\end{proof}

Now let $H_{1} = \rho_{1}(H), \; H_{2} = \rho_{2}(H) \leq \SU (n)$ be two realizations of $H$ as 
in Corollary~\ref{Cor:Reduction}.
Since $H_{1}, H_{2} \leq \SU (n)$ have the same dimension data with respect to the standard 
representation of $\SU(n)$ it follows from Frobenius' Reciprocity that $\quasiOne \sim \quasiTwo$.
Now, if we consider $\SU(n)$ with the bi-invariant metric, then it is clear that $\SU(n)$ acts on itself 
by isometries and that $\pi_{i} : \SU(n) \to \SU(n)/H_{i}$ (the projection mapping) is a Riemannian submersion
with totally geodesic fibers for $i = 1,2$. It then follows from Theorem~\ref{Thm:GenSuPe} that
the quotient spaces $\SU (n) /H_{1}$ and $\SU (n)/H_{2}$ are isospectral. 
From their construction as quotients of $\SU(n)$ it is clear from O'Neill's formula
\cite{ONeill} that these spaces have  non-negative sectional curvature.
We also note that it follows from the exact homotopy 
sequence of a weak fibration that these spaces are simply connected (see \cite[Chapter 4]{Switzer}).
We now turn our attention to the task of showing these spaces are locally non-isometric.

It is well known that simply-connected homogeneous spaces are isometric if and only if 
they are locally isometric. 
Consequently, it is enough to show that these spaces are non-isometric.
In \cite{On} the isometry groups of homogeneous spaces are studied and we see that for $i = 1,2$ the connected
component of the identity element, $\Isom(G/H_{i})^{0}$, is the locally direct product of $G$ and
$[N_{G}(H_{i})/H_{i}]^{0}$, which is denoted by $G\cdot [N_{G}(H_{i})/H_{i}]^{0}$. As $H_{i}$ and $G$ are connected and
$H_{i}$ is semi-simple it follows that $[N_{G}(H_{i})/H_{i}]^{0} \cong Z_{G}(H_{i})^{0}$ for each $i$. 
Hence, for $i = 1,2$ we have $$\Isom(G/H_{i})^{0} \cong G \cdot Z_{G}(H_{i})^{0}$$ and 
$$ \Isom(G/H_{i})^{0}_{\bar{e}_{i}} \cong H_{i} \cdot Z_{G}(H_{i})^{0},$$ 
where $\bar{e}_{i} = \pi_{i}(e)$ and $\pi_{i} : G \to G/H_{i}$ is the canonical projection for $i = 1,2$.

We now assume there is an isometry $f : (G/H_{1}, m_{1}) \to (G/H_{2}, m_{2})$. Without loss of generality 
we may assume that $f(\bar{e}_{1}) = \bar{e}_{2}$. The isometry $f$ then induces a Lie group isomorphism 
$\alpha : \Isom(G/H_{1})^{0} \to \Isom(G/H_{2})^{0}$ given by $\alpha(\Psi) = f \circ \Psi \circ f^{-1}$.
Since $\alpha$ must map simple factors to simple factors and $G = \SU(n)$ is a simple factor contained in 
neither $Z_{G}(H_{1})$ or $Z_{G}(H_{2})$ we conclude that $\alpha(G) = G$ and 
$\alpha(Z_{G}(H_{1})^{0}) = Z_{G}(H_{2})^{0}$. 
Also, since $\alpha (\Isom(G/H_{1})^{0}_{\bar{e}_{1}}) = \Isom(G/H_{2})^{0}_{\bar{e}_{2}}$ it follows that 
$\alpha(H_{1}) = H_{2}$. So we see that our isometry $f$ induces an automorphism $\alpha : G \to G$ such that
$\alpha(H_{1}) = H_{2}$, which is a  contradiction by Corollary~\ref{Cor:Reduction}. Hence, our spaces are not
isometric, and consequently they are locally non-isometric.

We may summarize our work thus far as follows.

\begin{theorem}\label{Thm:MyExamples}
There exists a connected, simply-connected semi-simple real Lie group $H$ which for infinitely many $n \in \N$ admits 
reducible faithful representations $\rho_{1}, \; \rho_{2} : H \to \SU(n)$, where $\rho_{1} \not\auto \rho_{2}$ and 
$H_{1} = \rho_{1}(H)$ and $H_{2} = \rho_{2}(H)$ are not conjugate by $\Aut(\SU(n))$.
If we equip $\SU(n)$ with a bi-invariant metric $m$, then the simply-connected, normal homogeneous spaces 
$(\SU(n)/H_{1}, m_{1})$ and $(\SU(n)/H_{2}, m_{2})$ are isospectral yet locally non-isometric.
\end{theorem}

\begin{remark} We offer the following comments.
\begin{enumerate}
\item It is clear that if one picks $\Gamma_{1}, \; \Gamma_{2} \leq \SU(n)$ discrete such that 
$[\pi^{G}_{\Gamma_{1}} : \rho] = [\pi^{G}_{\Gamma_{2}} : \rho]$ for all $\rho \in \widehat{\SU(n)}_{H_{1}} =\widehat{\SU(n)}_{H_{2}}$, then 
 $\Gamma_{1} \backslash \SU(n)/H_{1}$ and $\Gamma_{2} \backslash \SU(n)/H_{2}$ are isospectral yet locally non-isometric.
\item There is no $\SU(n)$-equivariant homeomorphism between $\SU(n)/H_{1}$ and $\SU(n)/H_{2}$. However, we cannot at this time 
determine whether the spaces are homeomorphic.
\item The smallest value of $n$ in Theorem~\ref{Thm:MyExamples} will be quite large, this follows 
from the comment on p. 393 of \cite{LP}. In fact we estimate that the dimension of the smallest resulting
homogeneous space is on the order of $10^{10}$.
\end{enumerate}
\end{remark}

\section{Isospectral Fiber Bundles}\label{Sec:FiberBundles}

In the previous section we saw that the study of dimension data can lead to examples of isospectral pairs 
which are quotients of compact Lie groups. In this section we will show that by considering dimension data 
we can also find isospectral pairs which  arise as quotients of Lie groups of non-compact type. Indeed we will 
establish the following result.

\begin{proposition}
Let $G$ be a semisimple Lie Group of non-compact type, $K \leq G$ a maximal compact subgroup and $\Gamma \leq G$ a 
co-compact lattice. Let $\rho : G \to \GL (V)$ be a finite dimensional faithful representation, so we may 
consider $G$ to be a closed linear group. Now suppose $H_{1}, H_{2} \leq K$ are closed, act freely on $G$ and 
have the same dimension data (with respect to $K$). It follows that $\Gamma \backslash G / H_{1}$ and 
$\Gamma \backslash G /H_{2}$ are isospectral on functions.
\end{proposition}

\noindent
The spaces $\Gamma \backslash G / H_{1}$ and $\Gamma \backslash G /H_{2}$ are fiber bundles over 
$\Gamma \backslash G / K$ with fibers  $(\Gamma \cap K) \backslash K / H_{1}$ and 
$(\Gamma \cap K) \backslash K /H_{2}$
respectively.

\begin{proof}
Endow $G$ with a metric which is left $G$-invariant and right $K$-invariant, hence when 
restricted to $K$ it is bi-invariant. 
Now select a co-compact lattice $\Gamma \leq G$. Then for any finite dimensional unitary representation 
$\sigma : K \to GL(V_{\sigma})$ of $K$ we may construct a locally homogeneous bundle 
$\pi: E_{\sigma} \to \Gamma \backslash G/K$. Indeed, let $K$ act on $(\Gamma \backslash G) \times V_{\sigma}$ by 
$k.(x,v) = (xk^{-1}, \sigma(k)v)$. Then we let $E_{\sigma} = ((\Gamma \backslash G) \times V_{\sigma})/K = 
\{ \overline{(x,v)} :
(x,v) \in (\Gamma \backslash G) \times V_{\sigma} \}$, where 
$\overline{(x,v)} = \{ (xk^{-1}, \sigma(k)v) : k \in K \}$. We let $\pi: E_{\sigma} \to \Gamma \backslash G/K$ be given by 
$\overline{(x,v)} \mapsto \overline{x}$. Then $\pi^{-1}(\overline{x}) = \{ \overline{(x,v)} : v \in V_{\sigma} \}$ is the 
fiber over $\overline{x} \in \Gamma \backslash G / K$.

We let $L^{2}(\Gamma \backslash G/K, E_{\sigma})$ denote the set of $L^{2}$-sections of $E_{\sigma}$. Then as a vector 
space $L^{2}(\Gamma \backslash G/K, E_{\sigma})$ is isomorphic to 
$$\widetilde{V}_{\sigma} \equiv \{ \widetilde{F} : \Gamma \backslash G \to V_{\sigma} : \widetilde{F} \mbox{ is }
L^{2},
\sigma(k)^{-1}\widetilde{F}(x) = 
\widetilde{F}(xk), \mbox{ for all } k \in K, \mbox{ a.e. } x \in \Gamma \backslash G \}.$$ 
This can be seen in the following manner.
Let $F \in  L^{2}(\Gamma \backslash G/K, E_{\sigma})$, then for all $\overline{x} \in \Gamma \backslash G/K$ 
we know $F(\overline{x}) \in \pi^{-1}(\overline{x})$. Hence, $F(\overline{x}) =
\overline{(x,\widetilde{F}(x))}$.  Now for $F$ to be well defined we must have for all $x \in \Gamma \backslash
G$ and for all $k \in K$ 
$(x, \widetilde{F}(x)) \sim (xk, \widetilde{F}(xk))$, but this occurs 
if and only if $\sigma(k)^{-1}\widetilde{F}(x) = \widetilde{F}(xk)$. So the correspondence is clear. On
$\widetilde{V}_{\sigma}$  we see that $K$ acts by $(k.\widetilde{F})(x) = \widetilde{F}(xk) =
\sigma(k)^{-1}\widetilde{F}(x)$.

We now recall that for any two measure spaces $(X, \mu)$ and $(Y, \nu)$ we have \linebreak 
$L^{2}(X\times Y) = L^{2}(X) \otimes L^{2}(Y)$ and if we have an action of a group $L$ on $X\times Y$ then 
\linebreak
$L^{2}(X\times_{L} Y) = (L^{2}(X) \otimes L^{2}(Y))^{L}$. Now, given that 
$\Gamma \backslash G = \Gamma \backslash G \times_{K} K$ and (by the Peter-Weyl Theorem) 
$L^{2}(K) = \oplus_{\sigma \in \widehat{K}} (\oplus_{i=1}^{\DIM \sigma} V_{\sigma})$ we see:
\begin{equation}
\begin{split}
L^{2}(\Gamma \backslash G) 	&= (L^{2}(\Gamma \backslash G) \otimes L^{2}(K))^{K} \\
				&= (L^{2}(\Gamma \backslash G) \otimes (\oplus_{\sigma \in \widehat{K}} (\oplus_{i =1}^{\DIM \sigma} V_{\sigma})))^{K} \\
				&= \oplus_{\sigma \in \widehat{K}} (\oplus_{i=1}^{\DIM \sigma} (L^{2}(\Gamma \backslash G) \otimes V_{\sigma})^{K}) \\
				&= \oplus_{\sigma \in \widehat{K}} (\oplus_{i=1}^{\DIM \sigma} L^{2}(\Gamma \backslash G /K, E_{\sigma}))\\
				&= \oplus_{\sigma \in \widehat{K}} (\oplus_{i=1}^{\DIM \sigma} \widetilde{V}_{\sigma}). \\
\end{split}
\end{equation}
Then for any $H \leq K$ we have 
$L^{2}(\Gamma \backslash G / H) = L^{2}(\Gamma \backslash G)^{H} = \oplus_{\sigma \in \widehat{K}}(\oplus_{i=1}^{\DIM \sigma}\widetilde{V}_{\sigma}^{H})$. 

In the case that $\sigma$ is a finite dimensional irreducible representation of $K$ we see that for 
$F \in \widetilde{V}_{\sigma}$ we 
know $\widetilde{F} \equiv 0$ or $V_{\sigma}= \mathcal{L}(Im (\widetilde{F}))$, the linear span of $Im(\widetilde{F})$.
Also, if 
$\widetilde{F} \in \widetilde{V}^{H}_{\sigma}$, $H \leq K$, then we see $Im \widetilde{F} \subset V^{H}_{\sigma}$. 
These facts imply that for $\sigma \in \widehat{K}$ and 
$\widetilde{F} \in \widetilde{V}_{\sigma}^{H} \backslash \{ 0\}$ we have 
$V_{\sigma} = \mathcal{L}(Im (\widetilde{F})) \subset V_{\sigma}^{H} \subset V_{\sigma}$, hence we conclude 
$V^{H}_{\sigma} = V_{\sigma}$ if and only if $\widetilde{V}^{H}_{\sigma} \neq 0$.
Therefore, for any $H \leq K$ we have $L^{2}(\Gamma \backslash G /H) = \oplus_{\{ \sigma \in \widehat{K}: \,
\Res^{K}_{H}(\sigma) = id \}}  (\oplus_{i=1}^{\DIM \sigma} \widetilde{V}_{\sigma}^{H})$.

We now recall that the bundle $E_{\sigma}$ admits a locally invariant connection $\nabla$, which is the push-forward 
of the invariant connection on the homogeneous bundle $\tilde{E}_{\sigma} = (G\times V_{\sigma})/K$. 
The connection $\nabla$ defines a quadratic 
form $D_{\sigma}$ on $C^{\infty}(\Gamma \backslash G /K, E_{\sigma})$ given by 
$$D_{\sigma}(f) = \int_{\Gamma \backslash G /K} \| \nabla f(x) \|^{2} dx. $$
The quadratic form $D_{\sigma}$ defines an elliptic operator $\Delta_{\sigma}$ on 
$L^{2}(\Gamma \backslash G/K, E_{\sigma})$ known as the Laplace operator. 
If $\sigma$ is irreducible, $\Delta_{\sigma}$ is equal to a shift of the restriction of the 
negative of the Casimir element of $G$ by a constant determined by $\sigma$. 
Now for any $H \leq K$ we see that $\Delta$ on $L^{2}(\Gamma \backslash G/H)$ is given by 
$\Delta = \oplus_{\{\sigma \in \widehat{K}: \; \Res^{K}_{H}(\sigma) = id \}} \Delta_{\sigma}$. 
It then follows from this and the above that if 
$H_{1}, H_{2} \leq K$ have the same dimension data (with respect to $K$), then 
$\Gamma \backslash G/H_{1}$ and $\Gamma \backslash G / H_{2}$ are isospectral.
\end{proof}

\begin{remark} 
If we let $H_{1}, \; H_{2} \leq K \equiv \SU(n)$ be as in Theorem~\ref{Thm:MyExamples}, $G = \SL_{n}(\C)$ 
and $\Gamma \leq G$ be co-compact, then we see that $\Gamma \backslash G /H_{1}$ and $\Gamma \backslash G / H_{2}$
are isospectral fiber bundles over $\Gamma \backslash G / K$ with isospectral fibers
$(\Gamma \cap K)\backslash K /H_{1}$ and $(\Gamma \cap K )\backslash K /H_{2}$.
\end{remark}


\section{Group Actions and a Theorem of von Neumann}\label{Sec:GroupActions}

We now conclude our paper by considering the spectra of group actions.
\medskip

Let $G$ be a locally compact group and $(X, \mu)$ a measure space where $X$ is 
a $G$-space and $\mu$ is a finite, $G$-invariant measure. We then obtain a representation of $G$ on $L^{2}(X,\mu)$ 
given by $(g \cdot f)(x) = f(g^{-1} \cdot x)$. The decomposition of $L^{2}(X,\mu)$ into $G$-irreducible
representations with their multiplicities taken into account is said to be the \textbf{spectrum} of the action of
$G$ on $X$.  If the decomposition of $L^{2}(X, \mu)$ into $G$-irreducibles is a countable direct sum of finite
dimensional irreducible representations we say that the spectrum of the action is {\bf discrete}.
Two $G$ actions are said to be \textbf{isospectral} if their spectra coincide.

A theorem of von Neumann states that two actions of a locally compact abelian group are measurably 
conjugate if their spectra are discrete and coincide \cite{vonNeumann}.
Spatzier considered the problem of spectral rigidity  of group actions in the case
of groups of non-compact type and obtained the following result.

\begin{theorem}[\cite{Sp}] \label{Spatzier}
Let $G$ be a non-compact almost simple connected real algebraic group whose complexification is one 
of the following types:
\begin{enumerate}
\item $A_{n}$ with $n \geq 26$
\item $B_{n}$ with $n \geq 27$
\item $B_{n}$ or $D_{n}$ with $n \geq 13$
\end{enumerate}
Then $G$ has properly ergodic actions which are isospectral yet not measurably conjugate.
\end{theorem}
\noindent
However, the spectra of these actions are necessarily non-discrete. In particular, 
if $G$ is of non-compact type and $(X, \mu)$ is a $G$-space, then the $G$-irreducibles which occur
in the decomposition of $L^{2}(X,\mu)$ are infinite dimensional. 
Using the examples constructed in Theorem~\ref{Thm:MyExamples} we can show that, in general, actions  with
\emph{discrete} spectra are not characterized up to measurable conjugacy by their spectra. Indeed, we obtain the
following result.

\begin{proposition} \label{Prop:IsoAction}
Let $G = \SU (n)$, $H_{1}$ and $H_{2}$ be as in Theorem~\ref{Thm:MyExamples}. 
Any dense subgroup $\Theta \leq G$ has actions on the measure spaces $(G/H_{1}, dx_{1})$ and $(G/H_{2}, dx_{2})$ 
with discrete spectra which are isospectral, but the actions are not measurably conjugate.
\end{proposition}

\begin{proof}
Let $\Theta \leq G$ be dense and let $G$ act on $G/H_{1}$ and $G/H_{2}$ in the usual way. 
We then get actions of $\Theta$ on $(G/H_{1}, dx_{1})$ and $(G/H_{2}, dx_{2})$. 
Let us suppose these $\Theta$-actions are measurably conjugate. That is, suppose there exists 
$F: G/H_{1} \to G/H_{2}$ a measurable isomorphism such that 
$F(\theta. x) = \theta.F(x)$ for all  $\theta \in \Theta$. 
When $\mathcal{A} = \{ f : G/H_{1} \rightarrow G/H_{2} \mbox{ measurable} \}$ is endowed with 
the topology of convergence in measure it is a standard Borel space and we have a natural (Borel) 
action of $G$ on $\mathcal{A}$ given by 
$(g.f)(x) = g.f(g^{-1}.x)$. It can be seen that for all $f \in \mathcal{A}$ $G_{f}$ (the stabilizer of $f$) 
is closed. Since $\Theta$ is dense and $\Theta \subset G_{F}$ we have $G_{F} = G$. So, $F$ is a $G$-map. 
The same can be said for $F^{-1}$.

Now there exists $L_{1}: G/H_{1} \to G/H_{2}$ continuous such that $F = L_{1}$ a.e. and 
there exists $L_{2}: G/H_{2} \rightarrow G/H_{1}$ continuous such that $F^{-1} = L_{2}$ a.e. Then 
$L_{2} \circ L_{1} = I$ a.e. and $L_{1} \circ L_{2} = I$ a.e., where $I$ denotes the identity.
From continuity we obtain equality everywhere.
Consequently, $L =L_{1}: G/H_{1} \rightarrow G/H_{2}$ is a homeomorphism which is also a $G$-map.

It is clear that $G_{\bar{e}_{1}} = H_{1}$, where $\bar{e}_{1} = eH_{1}$. Then, since $L$ is a 
$G$-map, we see $H_{1} \leq G_{f(\bar{e}_{1})} = H^{g}_{2}$ for some $g \in G$. From the fact that $L$ is 
also a homeomorphism we see $H^{g}_{2} \leq G_{\bar{e}_{1}} = H_{1}$. We have thus established 
that $H_{1}$ and $H_{2}$ are conjugate in $G$. However, by construction this is false. We are then led to 
conclude that the $\Theta$-actions are not measurably conjugate.

Since $\Theta \leq G$ is dense we know that the spectra of the $\Theta$-actions coincide with the spectra 
of the respective $G$-actions. By construction the $G$-actions on $(G/H_{1}, dx_{1})$ and $(G/H_{2}, dx_{2})$ 
have discrete spectra and are isospectral. Hence, the $\Theta$ actions have discrete spectra and are isospectral.
\end{proof}

From Proposition~\ref{Prop:IsoAction} it follows that there are arithmetic lattices which admit actions with
discrete spectra that are isospectral yet not measurably conjugate. Indeed, we recall the following
result. 

\begin{proposition}[Restriction of Scalars]\label{Prop:RestrictionScalars}
Let $F \subset \R$ be an algebraic number field with $d = [F:\Q] < \infty$ and let $\mathcal{O}$ be the
ring of integers in $F$. Now suppose $G \leq \SL (n, \R)$ is defined over $F$, and let 
$\sigma_{1} = id, \sigma_{2}, \ldots, \sigma_{d}$ denote the $d$ distinct (up to complex conjugation)
imbeddings of 
$F$ in $\C$. Then $G_{\mathcal{O}}$ imbeds as an arithmetic lattice in 
$$G^{\sigma_{1}} \times \cdots \times G^{\sigma_{d}}$$
via the natural embedding $g \stackrel{\phi}{\mapsto} (\sigma_{1}(g), \ldots , \sigma_{d}(g))$, where 
$G^{\sigma_{i}}$ denotes the Galois conjugate of $G$ by $\sigma_{i}$. 
Furthermore, if $G$ is simple, then $\phi (G_{\mathcal{O}})$ is irreducible. 
\end{proposition}

If we let  $F = \Q [\sqrt{2}]$, then $\SU (m,l)$ $(m+l \geq 2)$, the set of matrices in $\SL(m+l)$ which
preserve  the quadratic form $\sum_{i= 1}^{m} x_{i}^{2} - \sqrt{2}(\sum_{i= m+1}^{m+l} x_{i}^{2})$, is
defined over $F$. In this case the only non-trivial imbedding of $F$ inside $\C$ is given by 
$\sigma(x+\sqrt{2} y) = x - \sqrt{2}y$ and $\mathcal{O} = \Z[\sqrt{2}]$. 
From Proposition~\ref{Prop:RestrictionScalars} it follows that 
$\Gamma = \phi(\SU(m,l)_{\mathcal{O}})$ is an irreducible arithmetic lattice in 
$\SU(m,l) \times \SU(m+l)$. In the case where $\mathrm{min}(m,l) \geq 1$, we see that $\SU(m+l)$ is 
the maximal compact factor in $\SU(m,l) \times \SU(m+l)$ and hence it follows from the
irreduciblity of $\Gamma$ that $\pi(\Gamma)$ is \emph{dense} in 
$\SU(m+l)$, where $\pi : \SU(m,l)\times \SU(m+l) \to \SU(m+l)$ is the canonical projection.
If we now let $\SU(n)$, $H_{1}$, and $H_{2}$ be as in Proposition~\ref{Prop:IsoAction}, then we see that
there is an irreducible arithmetic lattice $\Gamma$ in $\SU(n-2, 2) \times \SU(n)$ which has actions
on 
$\SU(n)/H_{1}$ and $\SU(n)/H_{2}$ with discrete spectra that are isospectral yet not measurably
conjugate. 

\begin{remark}
After a more careful review of the literature we have recently learned that a counterexample 
to the spectral rigidity of group actions is contained in \cite{Mackey}. Mackey observes that if 
$(G, H_{1}, H_{2})$ is a triple of groups where $H_{1}$ and $H_{2}$ are non-conjugate, representation
equivalent subgroups of $G$, 
then the $G$-actions on $G/H_{1}$ and $G/H_{2}$ are isospectral, but not measurably
conjugate. He then gives an example of such a triple of groups
taken from \cite{Todd}, where $G = S_{16}$ is the permutation group on 16 elements and $H_{1}$
and $H_{2}$ are two order 16 subgroups.
\end{remark}

\addcontentsline{toc}{section}{References}

\newcommand{\etalchar}[1]{$^{#1}$}

\bigskip 
\noindent
Craig J. Sutton \\
Department of Mathematics \\
University of Pennsylvania \\
Philadelphia, PA 19104-6395 \\
USA \\
e-mail: cjsutton@math.upenn.edu     

\end{document}